\documentclass{amsart}
\usepackage[lmargin=3cm,rmargin=3cm,top=4cm,bottom=4cm]{geometry}
\newtheorem{theorem}{Theorem}[section]
\newtheorem{proposition}[theorem]{Proposition}
\newtheorem{lemma}[theorem]{Lemma}

\newtheorem{definition}[theorem]{Definition}
\newtheorem{remark}[theorem]{Remark}
\numberwithin{equation}{section}
\usepackage{caption} % For caption spacing
\usepackage{subcaption} % For sub-figures
\usepackage{graphicx}
\usepackage{pgfplots}
\usepackage[all]{nowidow}
\usepackage[utf8]{inputenc}
\usepackage{tikz}
\usepackage[inline]{enumitem}
\usepackage{blindtext}
\usepackage{amssymb}		
\usepackage{ulem}
\usepackage{amsfonts}

\usepackage[hidelinks,colorlinks=true,linkcolor=blue,citecolor=red, urlcolor=blue]{hyperref}
\numberwithin{equation}{section}
\pgfplotsset{compat=1.18}
\begin{document}
\title{Variational Proof of Conditional Expectation}

\author[Hugo Reyna-Castañeda]{Hugo Guadalupe Reyna-Castañeda}
\address{Departamento de Matemáticas, Facultad de Ciencias, Universidad Nacional Autónoma de México, Mexico City, Mexico}
\curraddr{}
\email{hugoreyna46@ciencias.unam.mx}
\thanks{}

%    author two information
\author{Mar\'ia de los \'Angeles Sandoval-Romero}
\address{Corresponding author: Departamento de Matemáticas, Facultad de Ciencias, Universidad Nacional Autónoma de México, Mexico City, Mexico}
\curraddr{}
\email{selegna@ciencias.unam.mx}
\thanks{}

\begin{abstract}
In this paper, we show that the conditional expectation of a random variable with finite second moment given a $\sigma$-algebra is the unique critical point of an energy functional in Hilbert space $L^2$. Then, we extend by density the result to every integrable random variable.
\end{abstract}
\maketitle

\setlength{\parskip}{0.25em} 
\section{Introduction} \label{Introduction}
\sloppy

Let $(\Omega,\mathcal{F},\mathbb{P})$ be a probability space. The concept of conditional expectation is fundamental to understanding the average behavior of a random variable, given certain conditions or additional information. It is especially used in the theory of martingales; see, for example, \cite[Chapter 10]{Gut}.

\begin{definition}\label{ConditionalExpectation}
Let $X:\Omega \to \mathbb{R}$ be an integrable random variable and $\mathcal{G}$ $\sigma$-algebra such that $\mathcal{G}\subset \mathcal{F}$. The conditional expectation of $X$ given $\mathcal{G}$ is the unique, almost surely, function $\mathbb{E}(X\,|\,\mathcal{G}):\Omega \to \mathbb{R}$ $\mathcal{G}$-measurable and integrable such that:
        $$
        \int_{B} X\,d\mathbb{P} = \int_{B} \mathbb{E}(X\,|\,\mathcal{G})\,d\mathbb{P}\,\,\,\,\,\,\,\,\,\,\,\forall\,B \in \mathcal{G}.
        $$
\end{definition}

We must note that in this text we consider two random variables identical if they are equal except on a set of probability zero. That is, $X=Y$ if $\mathbb{P}(X=Y)=\mathbb{P}(\{\omega \in \Omega \,:\, X(\omega)=Y(\Omega)\})=1$.

There are two common proofs of the existence and uniqueness of the conditional expectation of an integrable random variable $X$ given a $\sigma$-algebra $\mathcal{G}$ such that $\mathcal{G} \subset \mathcal{F}$. The first one uses the Radon-Nikodym theorem. Set $\mathbb{Q}(B):=\mathbb{E}(X\cdot 1_{B})$ for $B \in \mathcal{G}$. The finite signed measure $\mathbb{Q}$ is absolutely continuous with respect to the probability measure $\mathbb{P}$ restricted to $\mathcal{G}$. Hence $\mathbb{E}(X\,|\,\mathcal{G})$ is the Radon-Nikodym derivative of $\mathbb{Q}$ with respect to $\mathbb{P}$ (see, for example, \cite[Chapter V-Section 4]{Kolmogorov} or \cite[Chapter 10-Theorem 1.2]{Gut}). 

The second uses the concept of orthogonality in the Hilbert space $L^2(\Omega,\mathcal{F},P)$ with the inner product $\langle X, Y\rangle_{2}:=\mathbb{E}(X\cdot Y)$. Riesz-Fischer theorem \cite[Theorem 4.8]{Brezis} guarantees that the space $L^2(\Omega,\mathcal{G},\mathbb{P})$ is a closed linear subspace of $L^2(\Omega,\mathcal{F}, \mathbb{P})$. Thus, given $X \in L^2(\Omega,\mathcal{F},\mathbb{P})$, the conditional expectation $\mathbb{E}(X\,|\,\mathcal{G})$ is the orthogonal projection of $X$ onto $L^2(\Omega,\mathcal{G},\mathbb{P})$ (see \cite[Theorem 22.6]{Jacod} and \cite[Definition 23.5]{Jacod}). Then, it is shown that the space $L^2(\Omega,\mathcal{F},\mathbb{P})$ is dense in the space $L^1(\Omega,\mathcal{F},\mathbb{P})$ and through a limiting process, the concept of conditional expectation is extended to any integrable random variable (see \cite[Lemma 23.1]{Jacod} and \cite[Theorem 23.4]{Jacod}).

Our goal in this text is to present the existence of the conditional expectation of a random variable in $L^2(\Omega,\mathcal{F},\mathbb{P})$ as a variational problem; that is, we will show that the conditional expectation is the unique critical point of an energy functional $J:L^2(\Omega,\mathcal{G},\mathbb{P}) \to \mathbb{R}$. Then, using the same density argument, the result is obtained for every element of $L^1(\Omega,\mathcal{F},\mathbb{P})$. This is an extension of the work presented in \cite{ReynaSandoval}.

\section{Variational Formulation and Proof} \label{2}

Let $(\Omega,\mathcal{F},\mathbb{P})$ be a probability space and $\mathcal{G}$ be a $\sigma$-algebra such that $\mathcal{G} \subset \mathcal{F}$. The space $L^2(\Omega,\mathcal{F},\mathbb{P})$ with inner product $\langle X,Y\rangle_{2}:=\mathbb{E}(X\cdot Y)=\int_{\Omega} X\cdot Y\,d\mathbb{P}$ and the norm $\Vert X \Vert_{2}=\sqrt{\mathbb{E}(X^2)}$ is a Hilbert space over $\mathbb{R}$ \cite[Theorem 4.8]{Brezis}.

The function $\mathbb{P}$ restricted to $\mathcal{G}$ is a probability measure in $(\Omega,\mathcal{G})$, so the space $L^2(\Omega,\mathcal{G},\mathbb{P})$ is, in fact, a closed linear subspace of $L^2(\Omega,\mathcal{F},\mathbb{P})$ (see \cite[Theorem 4.8]{Brezis} and \cite[Definition 23.5]{Jacod}). Consequently, $L^2(\Omega,\mathcal{G},\mathbb{P})$ is a Hilbert space over $\mathbb{R}$. 

Given $X \in L^2(\Omega,\mathcal{F},\mathbb{P})$, our goal is to prove the existence of a unique function $\xi \in L^2(\Omega,\mathcal{G},\mathbb{P})$ such that:
$$
\int_{\Omega} X \cdot 1_{B}\,d\mathbb{P} = \int_{B} \xi\cdot 1_{B}\,d\mathbb{P}\,\,\,\,\,\,\,\,\,\,\,\forall\,B \in \mathcal{G}.
$$

It is clear that indicator function $1_{B} \in L^2(\Omega,\mathcal{G},\mathbb{P})$ for all $B \in \mathcal{G}$. In fact, the linear subspace generated by the indicator functions of the elements in $\mathcal{G}$ is dense in $L^2(\Omega,\mathcal{G},\mathbb{P})$.

\begin{lemma}\label{Simpledense}
    If $\mathcal{I}(\Omega,\mathcal{G},\mathbb{P}):=\{1_{B}\,:\,B \in \mathcal{G}\}$ then $\mbox{\rm span}\{\mathcal{I}(\Omega,\mathcal{G},\mathbb{P})\}$ is a dense subspace in $L^2(\Omega,\mathcal{G},\mathbb{P})$.
\end{lemma}

The elements of $\mathcal{S}(\Omega,\mathcal{G},\mathbb{P}):=\mbox{\rm span}\{\mathcal{I}(\Omega,\mathcal{G},\mathbb{P})\}$ are called simple random variables. The proof of the lemma \ref{Simpledense} is essentially based on the fact that every random variable is the point-wise limit of a sequence of simple random variables and on Lebesgue's dominated convergence theorem (see \cite[Lemma 1.1]{Gut} and \cite[Theorem 9.1]{Jacod}). With this fact, the following result is essential in the formulation of this text.

\begin{lemma} \label{ExpectationProduct}
    Let $X \in L^2(\Omega,\mathcal{G},\mathbb{P})$ and let $\mathbb{E}(X\,|\,\mathcal{G}) \in L^2(\Omega,\mathcal{G},\mathbb{P})$ be its conditional expectation. Then,
    $$
\int_{\Omega} X\cdot Y\,d\mathbb{P}=\int_{\Omega} \mathbb{E}(X\,|\,\mathcal{G})\cdot Y\,d\mathbb{P} \,\,\,\,\,\,\,\,\,\,\,\,\forall\,\,Y \in L^2(\Omega,\mathcal{G},\mathbb{P}).
    $$
\end{lemma}

\begin{proof}
    Let $Y \in L^2(\Omega,\mathcal{G},\mathbb{P})$. At first, we assume that $Y$ is a simple random variable, that is, $Y=\sum_{j=1}^{N} \alpha_j1_{B_j}$ with $\alpha_1,\ldots,\alpha_N \in \mathbb{R}$ and $B_1,\ldots,B_j \in \mathcal{G}$. By the linearity of the integral,
    $$
\int_{\Omega} \mathbb{E}(X\,|\,\mathcal{G})\cdot Y\,d\mathbb{P} =\sum_{j=1}^{N}\alpha_j\int_{\Omega} \mathbb{E}(X\,|\,\mathcal{G})\cdot 1_{B_j} \,d\mathbb{P}=\sum_{j=1}^{N}\alpha_j\int_{\Omega} X\cdot 1_{B_j}\,d\mathbb{P} = \int_{\Omega} X \cdot Y\,d\mathbb{P}.
    $$

     Now let $Y$ be an arbitrary function of $L^2(\Omega,\mathcal{G},\mathbb{P})$. By Lemma \ref{Simpledense} there exists $(Y_k)$ sequence of random variables in $\mathcal{S}(\Omega,\mathcal{G},\mathbb{P})$ such that $\lim_{k \to \infty}\Vert Y_k - Y \Vert_{2}=0$. Hence, there exist a subsequence $(Y_{k_\ell})$ and a function $Z\in L^2(\Omega,\mathcal{G},\mathbb{P})$ such that $\lim_{\ell \to \infty} Y_{k_\ell}=Y$ and $|Y_{k_\ell}| \leq Z$ almost surely (see \cite[Theorem 17.3]{Jacod} and \cite[Theorem 4.9]{Brezis}). Thus, $\lim_{\ell \to \infty} \mathbb{E}(X\,|\,\mathcal{G})\cdot Y_{k_{\ell}}= \mathbb{E}(X\,|\,\mathcal{G}) \cdot Y$ and $|\mathbb{E}(X\,|\,\mathcal{G}) \cdot Y_{k_\ell}| \leq |\mathbb{E}(X\,|\,\mathcal{G}) \cdot Z|$ almost surely with $|\mathbb{E}(X\,|\,\mathcal{G}) \cdot Z| \in L^1(\Omega,\mathcal{G},P)$ by the Hölder-Riesz inequality (see \cite[Theorem 9.3]{Jacod}). By Lebesgue’s dominated convergence theorem and the previous result
     $$
\lim_{\ell \to \infty}\int_{\Omega} X\cdot Y_{k_\ell}\,d\mathbb{P}= \lim_{\ell \to \infty}\int_{\Omega} \mathbb{E}(X\,|\,\mathcal{G})\cdot Y_{k_\ell} \,d\mathbb{P}= \int_{\Omega} \mathbb{E}(X\,|\,\mathcal{G})\cdot Y\,d\mathbb{P}
     $$

     In the same way, we can conclude that:
     $$
\lim_{\ell \to \infty}\int_{\Omega} X\cdot Y_{k_\ell}\,d\mathbb{P} = \int_{\Omega} X\cdot Y\,d\mathbb{P}
     $$
     therefore, necessarily:
     $$
\int_{\Omega} X\cdot Y\,d\mathbb{P}=\int_{\Omega} \mathbb{E}(X\,|\,\mathcal{G})\cdot Y\,d\mathbb{P}.
     $$
\end{proof}

\begin{remark} \label{InverseExpectationProblem}
    Now, given $X \in L^2(\Omega,\mathcal{F},\mathbb{P})$, if $\xi \in L^2(\Omega,\mathcal{G},\mathbb{P})$ is such that:
    $$
\int_{\Omega}X\cdot Y\,d\mathbb{P} = \int_{\Omega} \xi \cdot Y \,d\mathbb{P} \,\,\,\,\,\,\,\,\,\,\,\,\forall\,\,Y \in L^2(\Omega,\mathcal{G},\mathbb{P})
    $$
    then
    $$
\int_{\Omega}X\cdot 1_{B}\,d\mathbb{P} = \int_{\Omega} \xi \cdot 1_{B} \,d\mathbb{P} \,\,\,\,\,\,\,\,\,\,\,\,\forall\,\,B \in \mathcal{G}.
    $$
    and by uniqueness $\xi=\mathbb{E}(X\,|\mathcal{G})$.
\end{remark}

 By Lemma \ref{ExpectationProduct} and Remark \ref{InverseExpectationProblem}, we can present the existence of the conditional expectation as follows: Given $X \in L^2(\Omega,\mathcal{F},\mathbb{P})$, our goal is to prove the existence of a unique function $\xi \in L^2(\Omega,\mathcal{G},\mathbb{P})$ such that:
\begin{equation} \label{FProblem1}
    \int_{\Omega} X \cdot Y\,d\mathbb{P} = \int_{\Omega} \xi\cdot Y\,d\mathbb{P} =\langle \xi,Y\rangle_{2}\,\,\,\,\,\,\,\,\,\,\,\,\forall\,\,Y \in L^2(\Omega,\mathcal{G},\mathbb{P}).
\end{equation}

Since $X$ is given, we can write the left-hand side of the identity as a function in $L^2(\Omega,\mathcal{G},\mathbb{P})$. That is, define $T:L^2(\Omega,\mathcal{G},\mathbb{P}) \to \mathbb{R}$ by:
\begin{equation} \label{FProblem2}
    T(Y):=\int_{\Omega} X\cdot Y\,d\mathbb{P}.
\end{equation}

$T$ is clearly a linear function in $L^2(\Omega,\mathcal{G},P)$. Furthermore, by Hölder-Riesz inequality (see \cite[Theorem 9.3]{Jacod})
$$
|T(Y)| \leq \int_{\Omega} |X\cdot Y|\,d\mathbb{P} \leq \Vert X \Vert_{2}\Vert Y \Vert_{2}\,\,\,\,\,\,\,\,\,\,\,\,\forall\,\,Y \in L^2(\Omega,\mathcal{G},\mathbb{P})
$$
which implies that $T$ is continuous in $L^2(\Omega,\mathcal{G},P)$. Thus, we ask: Is there a unique function $\xi \in L^2(\Omega,\mathcal{G},P)$ such that $T(Y)=\langle \xi,Y\rangle_{2}$ for any $Y \in L^2(\Omega,\mathcal{G},P)$?

The previous question can be restated in the main result of this text. Theorem \ref{DirichletTheorem} states that the conditional expectation is the critical point of a functional in $L^2(\Omega,\mathcal{G},\mathbb{P})$ called the energy functional. This property is known as the Dirichlet principle (see \cite[Proposition 8.15]{Brezis}).

\begin{theorem} \label{DirichletTheorem}
Given $X \in L^2(\Omega,\mathcal{F},\mathbb{P})$, $\xi$ is the conditional expectation $\mathbb{E}(X\,|\,\mathcal{G})$ if and only if $\xi$ is the unique critical point of the functional $J:L^2(\Omega,\mathcal{G},\mathbb{P}) \to \mathbb{R}$ given by:
\begin{equation} \label{FunctionalEnergy}
    J(Y):=\frac{1}{2}\int_{\Omega} |Y|^2\,d\mathbb{P} - \int_{\Omega} X\cdot Y\,d\mathbb{P}.
\end{equation}
\end{theorem}

\begin{proof}
    $J$ is a $\mathcal{C}^2$ functional in $L^2(\Omega,\mathcal{G},\mathbb{P})$ such that (see Proposition \ref{DifLinear} and Proposition \ref{DifNorm})
    $$
J'_{Z}(Y)=\langle Z,Y\rangle_{2}-T(Y) \,\,\,\,\mbox{and}\,\,\,\,J_{Z}''(Y,W)=\langle Y,W\rangle_{2}
    $$
for all $W,Y,Z\in L^2(\Omega,\mathcal{G},\mathbb{P})$.

The proof is immediate from Lemma \ref{ExpectationProduct} and Remark \ref{InverseExpectationProblem}. In fact, $\xi$ is a critical point of $J$ if and only if $J'_{\xi}(Y)=0$ for all $Y \in L^2(\Omega,\mathcal{G},\mathbb{P})$ if and only if $T(Y)=\langle \xi,Y\rangle_{2}$ for all $Y \in L^2(\Omega,\mathcal{G},\mathbb{P})$.
\end{proof}

Thus, the proof of the concept of conditional expectation reduces to proving the existence of a unique critical point of $J$. The solution to problem (\ref{FProblem1}), which is the existence of a critical point of $J$, is given by the Fréchet-Riesz representation theorem in Hilbert spaces (see \cite[Theorem 5.5]{Brezis}). 

\begin{theorem}
    Given any $X \in L^2(\Omega,\mathcal{G},\mathbb{P})$, there is a unique critical point of functional $J:L^2(\Omega,\mathcal{G},\mathbb{P}) \to \mathbb{R}$ defined by:
    $$
J(Y):=\frac{1}{2}\int_{\Omega}|Y|^2\,d\mathbb{P}-\int_{\Omega} X\cdot Y \,d\mathbb{P}.
    $$
\end{theorem}

\begin{proof}
    Let $X \in L^2(\Omega,\mathcal{F},\mathbb{P})$. The function $T:L^2(\Omega,\mathcal{F},\mathbb{P}) \to \mathbb{R}$ given by 
    $$
T(Y)=\int_{\Omega} X \cdot Y\,d\mathbb{P}
    $$
is linear and continuous, so that, Fréchet-Riesz representation theorem ensures that there exists a unique element $\xi \in L^2(\Omega,\mathcal{G},\mathbb{P})$ such that:
$$
T(Y)=\langle \xi,Y\rangle_{2} \,\,\,\,\,\,\,\,\forall\,Y \in L^2(\Omega,\mathcal{G},\mathbb{P}).
$$

Consequently, $\xi$ is a critical point of $J$. Now, $J:L^2(\Omega,\mathcal{G},\mathbb{P})\to \mathbb{R}$ is of class $\mathcal{C}^2$ with $J_{\xi}''(Y,Z)=\langle Y,Z\rangle_{2}$ for all $Y,Z \in L^2(\Omega,\mathcal{G},\mathbb{P})$. Then, $J_{\xi}''(Y,Y)=\langle Y,Y\rangle_{2}=\Vert Y \Vert_{2}^{2} >0$ for all $Y \in L^2(\Omega,\mathcal{F},\mathbb{P})$ with $Y \neq 0$ which implies that $\xi$ is the unique minimum of $J$ (see Theorem \ref{SufMin}).
\end{proof}

We therefore conclude that the conditional expectation of $X$ given $\mathcal{G}$ is the unique minimum of the energy functional $J$ defined in (\ref{FunctionalEnergy}).

\appendix
\section{}

In this brief appendix, we will gather concepts on the differentiability of the energy functional, which are generally developed for an arbitrary Hilbert space over $\mathbb{R}$. For more details, we suggest consulting \cite{Drabek} and \cite{Gelfand}.

Let $H=(H,\langle \cdot,\cdot\rangle,\Vert \cdot \Vert_H)$ be a Hilbert space over $\mathbb{R}$, $H^\ast:=\{L:H\to\mathbb{R}\,:\,L\,\,\mbox{is linear and continuos}\}$ a dual space, and $\mathcal{B}(H \times H):=\{L:H \times H \to \mathbb{R}\,:\,L\,\,\mbox{is bilinear and continuos}\}$.

\begin{definition} \label{FirstDerivate}
    Let $F: H \to \mathbb{R}$ be a function. 
    \begin{itemize}
        \item[(a)]  If, for $u \in H$ and $v \in H$, the limit
$$
\lim_{t \to 0} \frac{F(u+tv)-F(u)}{t}
$$
exists, then its value is called the derivative of $F$ at the point $u$ and in the direction $v$.

\item[(b)] $F:H \to \mathbb{R}$ is Gateaux-differentiable at the point $u_{0} \in H$ if, for every $v \in H$, the derivative of $F$ at the point $u$ and in the direction $v$ exists and the function $F'_{u}:H \to \mathbb{R}$ given by
$$
F'_{u}(v):=\lim_{t \to 0} \frac{F(u+tv)-F(u)}{t}
$$
is an element on $H^{\ast}$.

\item[(c)] $F:H \to \mathbb{R}$ is Gateaux-differentiable on $H$ if it is Gateaux-differentiable at every point $u \in H$. The function given by
$$
F':H\to H^{\ast},\,\,\,\,\,\,\,\,u \mapsto F'_{u},
$$
is called the Gateaux derivative of $F$. 

\item[(d)] $F:H \to \mathbb{R}$ is of class $\mathcal{C}^1$ in $H$ if and only if $F$ is Gateaux-differentiable on $H$ and $F':H \to H^{\ast}$ is continuous (see \cite[Proposition 3.2.15]{Drabek}). 

\item[(e)] $u_0$ is a critical point of $F$ in $H$ if $F'_{u_0}(v)=0$ for all $v \in H$. 
\end{itemize}
\end{definition}

\begin{definition} \label{SecondDerivate}
    Let $F:H \to\mathbb{R}$ be a function of class $\mathcal{C}^1$.
\begin{itemize}
    \item[(a)] $F$ has a second Gateaux derivative at $u_0 \in H$ if there exists a continuous bilinear form $L :H  \times H \to \mathbb{R}$ such that for every $v,w \in V$
    $$
     \lim_{t \to 0}\frac{ F'_{u_0+tv}(w) - F'_{u_{0}}(v) - L(tv,w) }{t}=0.
    $$
     
     We denote a second Gateaux derivative of $F$ in $u_0$ by $L:=F''_{u_{0}}$. Hence, the second Gateaux derivative on $u_{0} \in H$ is given by (see \cite[Remark 3.2.29]{Drabek})
$$
F''_{u_{0}}(v,w):=\lim_{t \to 0} \frac{F'_{u_{0}+tv}(w) - F'_{u_{0}}(w)}{t}.
$$

     \item[(b)] $F:H \to \mathbb{R}$ is twice Gateaux-differentiable on $H$ if it has a second Gateaux derivative at every point $u \in H$. The function given by
$$
F'':H \to \mathcal{B}(H \times H),\,\,\,\,\,\,\,\,u \mapsto F''_{u},
$$
is called the second Gateaux derivative of $F$.

\item[(c)] $F$ is of class $\mathcal{C}^2$ in $H$ if and only if $F$ has a continuous second Gateaux derivative on $H$.
\end{itemize}
\end{definition}

According to the established objective, the following result provides a sufficient condition for finding minimum points of certain functionals \cite[Theorem 2]{Gelfand}.

\begin{theorem} \label{SufMin}
Let $F:H \to \mathbb{R}$ be a function of class $\mathcal{C}^{2}$. $F$ has a minimum point $u \in H$ if $F'_{u}(v)=0$ for all $v \in H$ and $F''_{u}(v,v) > 0$ for all $v \in H \smallsetminus \{0_H\}$.
\end{theorem}

The energy functional $J$ in (\ref{FunctionalEnergy}) is defined in terms of the squared norm of the Hilbert space $L^2(\Omega,\mathcal{G},\mathbb{P})$ and the continuous linear function $T:L^2(\Omega,\mathcal{G},\mathbb{P}) \to \mathbb{R}$. In the following propositions, we show that the energy functional is of class $\mathcal{C}^2$ and calculate its derivatives.

\begin{proposition}\label{DifLinear}
If $L \in H^{\ast}$ then $L$ is of class $\mathcal{C}^{2}$ in $H$, $L'=L$ and $L''=0$.
\end{proposition}

\begin{proof}
Let $u \in H$. We have $L(u+v)-L(u)-L(v)=0$ for all $v \in H$. Thus, $L'_{u}=L$ for all $u \in H$. Given that $L'_{u+v}=L$ for all $u,v \in H$ then $L'$ is a constant function and, thus, continuous. Consequently, $L''=0$ and $L$ are of class $\mathcal{C}^2$.
\end{proof}

\begin{proposition} \label{DifNorm}
    The function $F:H \to \mathbb{R}$ given by $F(u)=\frac{1}2\Vert u \Vert_{H}^{2}$ is of class $\mathcal{C}^2$, $F'_{u}(v)=\langle u,v \rangle$ and $F''_{u}(v,w)=\langle v,w \rangle$ for all $v,w \in H$.
\end{proposition}

\begin{proof}
Let $u \in H$. For every $v \in H$ and $t \in \mathbb{R}$ we have that:
$$
\begin{aligned}
F(u+tv)-F(u) &=\frac{1}{2}\left( \Vert u \Vert_{H}^2 - 2t\langle u,v\rangle + t^2\Vert v \Vert_{H}^2  - \Vert u \Vert_{H}^2 \right)= t\langle u,v\rangle + \frac{t^2}{2}\Vert v \Vert_{H}^2.
\end{aligned}
$$

Hence, 
$$
\lim_{t \to 0} \frac{F(u+tv)-F(u)}{t} = \lim_{t \to 0} \left( \langle u,v \rangle - \frac{t}{2}\Vert v \Vert_H^2 \right) = \langle u,v\rangle \,\,\,\,\,\,\,\forall\,v \in H.
$$

That is, $F'_{u}(v)=\langle u,v \rangle$ for all $v \in H$ and it is clear that $F'_{u} \in H^{\ast}$ by Fréchet-Riesz representation theorem (see \cite[Theorem 5.5]{Brezis}). Let $(u_j)$ be a sequence on $H$ such that $\Vert u_j - u \Vert_H \to 0$ as $j \to \infty$. By Cauchy–Schwarz inequality,
$$
|F'_{u_j}(v)-F'_{u}(v)|=|\langle u_j - u, v \rangle| \leq \Vert u_j - u \Vert_H \Vert v \Vert_H \,\,\,\,\,\forall\, v \in H.
$$

Thus, $F'$ is continuous and $F$ is of class $\mathcal{C}^1$ on $H$.

Now, for every $v,w \in H$ and $t \in \mathbb{R}$, $F'_{u+tv}(w)-F'_{u}(w)=\langle u+tv,w\rangle - \langle u,w \rangle = t\langle v,w \rangle$ and consequently:
$$
\lim_{t \to 0} \frac{F'_{u+tv}(w)-F'_{u}(w)}{t} = \langle v,w \rangle\,\,\,\,\,\,\,\,\,\,\,\,\,\,\forall v,w \in H.
$$

Therefore, $F''_{u}(v,w)=\langle v,w \rangle$ for all $v,w \in H$ and is a bilinear continuous form. This shows that $F$ is of class $\mathcal{C}^2$.
\end{proof}

\bibliographystyle{amsplain}

\begin{thebibliography}{10}

\bibitem{Gut} {A. Gut}, \textit{Probability: A Graduate Course. Second Edition}. Springer,  New York, 2013.

\bibitem{Kolmogorov} {A.N. Kolmogorov}, \textit{Foundations of the Theory of Probability}. Chelsea Publishing Company, New York, 1956.

\bibitem{Brezis} {H. Brezis}, \textit{Functional Analysis, Sobolev Spaces and Partial Differential Equations}. Springer,  New York, 2010.

\bibitem{ReynaSandoval} {H. Reyna and Ma. de los A. Sandoval}, \textit{A variational method to calculate probabilities}. arXiv: \href{https://arxiv.org/abs/2503.16727}{2503.16727}, 2025.

\bibitem{Gelfand} {I.M. Gelfand and S.V. Fomin}, \textit{Calculus of Variations}. Prentice-Hall, Inc. Englewood Cliffs, New Jersey, 1963.

\bibitem{Jacod} {J. Jacod and P. Protter}, \textit{Probability Essentials}. Springer, Berlin, 2004.

\bibitem{Drabek} {P. Drábek and J. Milota}, \textit{Methods of Nonlinear Analysis. Aplications to Differential Equations}. Birkh\"auser, 2007.


\end{thebibliography}

\end{document}